
\input amssym 
\magnification = \magstep1

\centerline{\bf ON EIGENVALUES OF THE KERNEL ${1\over 2} +\lfloor {1\over xy}\rfloor - {1\over xy}\,$ ($0<x,y\leq 1$)}

\bigskip

\centerline{par Nigel Watt}

\vskip 7mm

\noindent{\bf R\'{e}sum\'{e}.}\quad 
Nous montrons que le noyau $K(x,y) 
= {1\over 2} + \lfloor {1\over x y}\rfloor - {1\over x y}\,$ ($0<x,y\leq 1$)   
a un nombre infini de valeurs propres positives 
et un nombre infini de valeurs propres n\'{e}gatives. Notre int\'{e}r\^{e}t pour 
ce noyau est motiv\'{e} par l'apparition de la forme quadratique 
$\sum_{m,n\leq N} K\bigl( {m\over N} , {n\over N}\bigr) \mu(m)\mu(n)$ 
dans une identit\'{e} impliquant la fonction de Mertens. 

\medskip

\noindent{\bf Abstract.}\quad We show that the kernel $K(x,y) 
= {1\over 2} + \lfloor {1\over x y}\rfloor - {1\over x y}\,$ ($0<x,y\leq 1$) 
has infinitely many positive eigenvalues and infinitely many negative 
eigenvalues. Our interest in this kernel is motivated by the appearance 
of the quadratic form $\sum_{m,n\leq N} K\bigl( {m\over N} , {n\over N}\bigr) \mu(m)\mu(n)$ 
in an identity involving the Mertens function. 
\footnote{}{2010 {\it Mathematics Subject Classification.} 11A25, 45C05, 11A07.}
\footnote{}{{\it Mots-clefs.} Mertens function, eigenvalue, symmetric kernel.}

\vskip 7mm

\noindent{\bf \S 1. Introduction.}

\medskip 

For  $0<x,y\leq 1$, put 
$$K(x,y)={1\over 2} - \left\{ {1\over xy}\right\}\;,$$
where $\{ t\}\in[0,1)$ denotes the fractional part of $t\in{\Bbb R}$ (i.e. 
$\{ t\} = t -\lfloor t\rfloor$, where $\lfloor t\rfloor = \max\{ m\in{\Bbb Z} \,:\, m\leq t\}$). 
When $0\leq x,y\leq 1$ and $xy=0$, put $K(x,y)=0$. The function $K$ thus defined on 
$[0,1]\times[0,1]$ is (in the terminology of [3]) a symmetric, non-null $L_2$-kernel. 
It is shown in [3, Section~3.8] that every such kernel has at least one eigenvalue~$\lambda$. 
That is, there exists a number $\lambda\neq 0$, and an associated `eigenfunction' 
$\phi(x)$ (with $\infty > \int_0^1 |\phi (x)|^2 dx > 0$), satisfying
$$\phi(x) = \lambda\int_0^1 K(x,y)\phi(y) dy \eqno(1.1)$$ 
almost everywhere, with respect to the Lebesgue measure, in $[0,1]$. 
\par 
Since $K$ is symmetric (i.e. satisfies $K(x,y)=K(y,x)$), all eigenvalues of 
$K$ are real, and so there is no essential loss of generality in 
considering just those eigenfunctions of $K$ that are real valued 
(i.e. at least one of the pair of real functions ${\rm Re}(\phi(x))$, ${\rm Im}(\phi(x))$ 
may be substituted for $\phi(x)$ in (1.1)). 
\par 
In this paper we have the option of working only with eigenfunctions 
$\phi : [0,1]\rightarrow {\Bbb R}$ that satisfy (1.1) for all $x\in [0,1]$ 
(when $\lambda\in{\Bbb R}\backslash\{ 0\}$ is the appropriate eigenvalue).  
Choosing to do so would not be overly restrictive, for if $\phi(x)$ is any   
eigenfunction of $K$, with associated eigenvalue $\lambda$, then (1.1) holds almost everywhere in $[0,1]$,  
and the term  $\lambda \int_0^1 K(x,y) \phi(y) dy$ (occurring in (1.1))    
is an eigenfunction of $K$ that has the required property   
(the last part of this following, via the Cauchy-Schwarz inequality, from the observation that,  
for all $x\in [0,1]$, the integral $\int_0^1 (K(x,y))^2 dy$ exists, and is finite). 
Although this is an option that is 
of no consequence with regard to the proof of our main result (Theorem~1, below), 
we shall find it helpful when discussing certain incidental matters. 
\par 
By the general theory set out in [3, Section~3.8],
the set 
$${\cal S}(K)=\{ \lambda\,:\, \lambda\ {\rm is\ an\ eigenvalue\ of}\ K\}$$ 
is countable (in the sense that does not preclude its being finite). 
It is not hard to see that ${\cal S}(K)$ cannot be a finite set  
(for a sketch of a proof of this, see our Remark~(2), following Lemma~2.1 below).  
This paper is devoted to proving the following stronger result. 

\medskip 

\noindent{\bf THEOREM~1.}\ {\it Each of sets ${\cal S}(K)\cap (-\infty ,0)$ and ${\cal S}(K)\cap (0, \infty )$ is infinite.} 

\medskip 

We prove this theorem in Section~3, after some necessary preliminaries. 

Our particular interest in the kernel $K(x,y)$ is motivated by 
a connection with the M\"obius function $\mu(n)$ and its associated 
summatory function $M(x)=\sum_{m\leq x} \mu(n)$ (known as the 
Mertens function). This connection is apparent in our recent joint work [1] with Huxley, 
where it is (in effect) noted that for each positive integer $N$ one has  
$${M\left( N^2\right) \over N^2} 
+ \sum_{m=1}^N \sum_{n=1}^N K\left( {m\over N} , {n\over N}\right) {\mu(m) \mu(n) \over N^2} 
= {M(N) \left( M(N) + 4\right) \over 2 N^2} - \left( \sum_{m=1}^N {\mu(m)\over m}\right)^{\!\!2} 
\eqno(1.2)$$
(this following directly from [1, Equations~(3)-(5) and~(37)]). 
In the work [4] (to appear) it is shown that the sum over $m$ and $n$ on 
the left-hand side of (1.2) may 
be approximated, reasonably well, by sums involving the numbers $\mu(1),\ldots ,\mu(N)$, 
certain of the (smaller) eigenvalues of $K$ 
and the values that the corresponding eigenfunctions have at the points 
$x=\textstyle{m\over N}\,$ ($m=1,\ldots ,N$). 
\par
The author is indebted to the anonymous referee for pointing out that F. Mertens himself 
showed (in 1897) that, for all positive integers $n$, one has    
$$M(n) = 2 M\left( \sqrt{n}\right) - 
\mathop{\sum\sum}_{r,s\leq \sqrt{n}} \mu(r) \mu(s) \left\lfloor {n\over rs}\right\rfloor .$$ 
The proof of this appeared in [2, Section~3]. 
This result of Mertens contains the `principal case' of [1, Equations~(3)-(5)], 
from which we have derived the equation (1.2),  
and is equivalent to that subcase of the `principal case' of [1, Theorem~1] in which 
one has $d=2$ and $N_1=N_2=\lfloor \sqrt{K}\rfloor$ (with $K=n$). 

\bigskip

\noindent{\bf Acknowledgement.}\quad The author wishes to thank the anonymous referee 
for pointing out the relevance of the work [2] of Mertens, and for 
several other comments that have helped to improve this paper. 

\bigskip

\noindent{\bf \S 2. Notation and some Hilbert-Schmidt Theory.}

\medskip

\par  
We denote by $L^2 ([0,1])$ the semimetric space of 
functions $f : [0,1] \rightarrow {\Bbb R}$ that are measurable (in the sense of Lebesgue) 
and satisfy the condition $\int_0^1 (f(x))^2 dx < \infty$.  
\par 
For each eigenvalue $\lambda\in{\cal S}(K)$, 
the corresponding set of eigenfunctions (satisfying (1.1) for all $x\in [0,1]$) spans a finite 
dimensional subspace of the space $L^2 ([0,1])$: we follow [3] in referring to the dimension, $r_{\lambda}\,$ (say), 
of this subspace as the `index' of~$\lambda$.  
We put $\omega = \sum_{\lambda\in\cal S} \, r_{\lambda}$, so that $\omega\in{\Bbb N}$ if ${\cal S}(K)$ is 
finite, while $\omega =\infty$ otherwise. 
Since $K$ is symmetric, eigenfunctions $\phi,\psi$ (say) corresponding 
to two different eigenvalues $\lambda,\mu\in{\cal S}(K)$ are necessarily orthogonal with 
respect to the (semi-definite) inner product 
$$\langle f , g\rangle = \int_0^1 f (x) g (x) dx 
\qquad\hbox{($f,g\in L^2 ([0,1])$).}\eqno(2.1)$$
See [3, Sections~3.1 and~2.3] regarding this matter. 
In [3, Section~3.8] it is shown that there exists a system 
$\phi_j\,$ ($j\in{\Bbb N}$ and $j\leq\omega$) of eigenfunctions of $K$ that is 
orthonormal, so that one has  
$$\langle \phi_j , \phi_k\rangle :=\int_0^1 \phi_j(x) \phi_k(x) dx 
=\cases{1 &if $j=k$, \cr 0 &otherwise,\cr}$$
whenever $j,k\in{\Bbb N}$ satisfy $j,k\leq\omega$, 
and that is (at the same time) maximal, so that the corresponding sequence 
$\lambda_j\,$  ($j\in{\Bbb N}$ and $j\leq\omega$) of eigenvalues of $K$ is such that  
one has $|\{ j\in{\Bbb N} \,:\, j\leq\omega\ {\rm and}\ \lambda_j = \lambda\}|=r_{\lambda}$ 
for  all $\lambda\in{\cal S}(K)$. 
\par 
By [3, (8),~Section~3.10], we have 
$$\sum_{\scriptstyle j\in{\Bbb N}\atop\scriptstyle j\leq\omega} {1\over \lambda_j^2} 
=\int_0^1 \int_0^1 \left( K(x,y)\right)^2 dx dy <{1\over 4}\;.\eqno(2.2)$$
It follows that either ${\cal S}(K)$ is finite, or else one has 
$|\lambda_j|\rightarrow \infty$ as $j\rightarrow\infty$. Therefore, as in [3], 
we may assume that the $\phi_j$'s (and associated $\lambda_j$'s) 
are numbered in such a way that 
the absolute values of the associated eigenvalues form a sequence, 
$|\lambda_j|\,$ ($j\in{\Bbb N}$ and $j\leq\omega$), that is monotonic increasing. 
\par
We now develop some notation in which there is a clear 
distinction between positive and negative eigenvalues (and between 
the corresponding eigenfunctions). 
Let $\omega^{+}$ (resp. $\omega^{-}$) be the number of positive (resp. negative) 
terms in the sequence $\lambda_j\,$ ($j\in{\Bbb N}$ and $j\leq\omega$), 
so that $\omega^{+},\omega^{-}\in{\Bbb N}\cup\{0,\infty\}$ 
and $\omega^{+} + \omega^{-}=\omega$.
If all of the negative eigenvalues are removed from the sequence 
$\lambda_j\,$ ($j\in{\Bbb N}$ and $j\leq\omega$) 
then what remains is some monotonic increasing subsequence 
$\lambda_{m_k}\,$ ($k\in{\Bbb N}$ and $k\leq\omega^{+}$) in which each  
positive eigenvalue of $K$ appears.  
If one instead removes the positive eigenvalues then what remains is some   
monotonic decreasing subsequence $\lambda_{n_k}\,$ ($k\in{\Bbb N}$ and $k\leq\omega^{-}$) 
in which each negative eigenvalue of $K$ appears.
For $k\in{\Bbb N}$ satisfying $k\leq\omega^{+}\,$ (resp. $k\leq\omega^{-}$) 
we put $\lambda_k^{+}=\lambda_{m_k}$ and $\phi_k^{+}=\phi_{m_k}\,$ 
(resp. $\lambda_k^{-}=\lambda_{n_k}$ and $\phi_k^{-}=\phi_{n_k}$): note this 
has the consequence that (1.1) holds when $\lambda$ and $\phi$ are $\lambda_k^{+}$ and  $\phi_k^{+}\,$ 
(resp. $\lambda_k^{-}$ and  $\phi_k^{-}$), respectively. 
As every eigenvalue of $K$ is real and non-zero (and so either positive 
or negative), it is clear that the sets 
$\{\phi_j : j\in{\Bbb N}\}$ 
and $\{ \phi_k^{+} : k\in{\Bbb N}\ {\rm and}\ k\leq\omega^{+}\}
\cup \{ \phi_k^{-} : k\in{\Bbb N}\ {\rm and}\ k\leq\omega^{-}\}$ 
are equal, and so we know (in particular) that the elements of the latter set  
form an orthonormal system.  

\bigskip 

\noindent{\bf Lemma~2.1.}\quad 
{\it Let $\phi,\psi\in L^2 ([0,1])$ and put 
$$J=J(\phi,\psi)=\int_0^1 \int_0^1 K(x,y)\phi(x)\psi(y) dx dy\;.$$
Then 
$$J=\sum_{\scriptstyle j\in{\Bbb N}\atop\scriptstyle j\leq\omega} 
{\langle \phi, \phi_j\rangle \langle \phi_j, \psi\rangle\over \lambda_j}
=\sum_{\scriptstyle k\in{\Bbb N}\atop\scriptstyle\ k\leq\omega^{+}}
\!\!\!{\langle \phi, \phi_k^{+}\rangle \langle \phi_k^{+}, \psi\rangle\over \lambda_k^{+}} 
+ \sum_{\scriptstyle k\in{\Bbb N}\atop\scriptstyle\ k\leq\omega^{-}} 
\!\!\!{\langle \phi, \phi_k^{-}\rangle \langle \phi_k^{-}, \psi\rangle\over \lambda_k^{-}}\;.$$ 
}

\smallskip
\noindent{\bf Proof.}
\quad See [3, Section~3.11], where this result is proved in greater generality 
(i.e. for an arbitrary symmetric $L_2$-kernel) by applying a theorem of 
Hilbert and Schmidt (for which see [3, Section~3.10]) $\blacksquare$ 

\medskip 

\noindent{\bf Remarks:}

\smallskip

\noindent{\bf (1)}\ The proof supplied in [3] shows that 
each sum over $j$, or $k$, in Lemma~2.1 is absolutely convergent 
(when not finite or empty).   
This may also be deduced directly from Bessel's inequality,  
since, for $j\in{\Bbb N}$ with $j\leq\omega$, one has $|\lambda_j| > |\lambda_1| >0$  
and so, by the inequality of arithmetic  and geometric means, 
$|\langle \phi, \phi_j\rangle \langle \phi_j, \psi\rangle / \lambda_j| 
\leq 
\left( |\langle \phi, \phi_j\rangle|^2 + |\langle\phi_j, \psi\rangle|^2\right) 
/ (2 |\lambda_1|)$.  

\smallskip

\noindent{\bf (2)}\ Given what was noted in the third paragraph of Section~1, 
the eigenfunctions $\phi_j$ ($j\in{\Bbb N}$ and $j\leq\omega$) can be chosen in 
such a way that each has 
$\phi_j(x) = \lambda_j \int_0^1 K(x,y)\phi_j (y) dy$ for 
all $x\in [0,1]$. Assume that such a choice has been made.  
One can show that it follows that each $\phi_j$ is continuous on the interval $(0,1]$: we leave the proof 
of this as an exercise for the reader, and remark that 
the eigenfunctions in question can also be shown to be continuous at the 
point $x=0$ (this last fact, however, is not relevant 
to our main concern here). 
Therefore if $\omega\neq\infty\,$ (so that $\omega\in{\Bbb N}$)  
then,  by applying Lemma~2.1 for functions 
$\phi$ and $\psi$ that are positive valued and supported in intervals 
$[x-\varepsilon , x]$ and $[y-\varepsilon , y]$ (respectively), 
we find (letting $\varepsilon \rightarrow 0+$, and making use of the right-continuity 
of the real function $t\mapsto \{ t\}$) that one has 
$K(x,y)=\sum_{j=1}^{\omega} \phi_j(x)\phi_j(y) / \lambda_j$ at all points 
$(x,y)\in (0,1]\times (0,1]$.  
This last identity would imply that $K(x,y)$ is continuous on $(0,1]\times (0,1]$, 
whereas one has (for example)  
$\lim_{x\rightarrow \left( {1\over 2}\right) +} K(x,x) = -\textstyle{1\over 2}$, 
but $K(\textstyle{1\over 2} , \textstyle{1\over 2})=\textstyle{1\over 2}$.  
By this {\it reductio ad absurdum} we may conclude that $\omega =\infty$. 

\bigskip

\noindent{\bf \S 3. Negative (resp. positive) eigenvalues of $K$.}

\medskip 

In this section we prove Theorem~1 by showing that $\omega^{+}=\omega^{-}=\infty$. 
In doing so we shall make use of both 
Lemma~2.1 and the following purely number-theoretic result.

\bigskip 

\noindent{\bf Lemma~3.1.}\quad 
{\it 
Let $u\in\{ 1, -1\}$, let $Q\in [5,\infty)$ and let $N$ be a non-negative integer. 
Then there exist $N+1$ distinct primes 
$p_1,\ldots, p_{N+1}$, all greater than $Q$, and an integer 
$n$ satisfying 
$$3 n^2 \equiv m_j \pmod{p_j^2}\qquad \hbox{($j=1,\ldots ,N+1$),}\eqno(3.1)$$
where $m_j$ denotes the least positive integer satisfying both  
$$2 m_j\equiv 3\pmod{p_j}\qquad{\rm and}\qquad m_j\equiv u\pmod{3}\;.\eqno(3.2)$$
One has here 
$$0<m_j<3p_j\qquad\hbox{($j=1,\ldots ,N+1$),}\eqno(3.3)$$
and the integer $n$ may be chosen so as to satisfy 
$$P^2 < n < 2 P^2,\eqno(3.4)$$
where $P$ is the product of the primes $p_1,\ldots ,p_{N+1}$. 
}

\medskip

\noindent{\bf Proof.}\quad 
It is a corollary of Dirichlet's theorem on primes in arithmetic progressions that the 
set $\{ p : p\ {\rm is\ prime},\ p\equiv \pm 1\pmod{8}\ {\rm and}\ p>Q\}$ is  
infinite: we take $p_1,\ldots , p_{N+1}$ to be 
any $N+1$ distinct elements of this set. As $Q > 3$, and as 
distinct positive primes are coprime to one another, 
we have $(p_j,2)=(p_j,3)=1$ ($j=1,\ldots ,N+1$), and 
$(p_j,p_k)=1$ ($1\leq j<k\leq N+1$). It therefore follows by the Chinese Remainder Theorem 
that, for $j=1,\ldots,N+1$, the simultaneous congruences in (3.2) are soluble (for the 
integer $m_j$): since $(u,3)=(\pm 1, 3)=1=(3,p_j)$, the set of all integer solutions 
of these congruences 
is one of the residue classes $\bmod\ 3p_j$ that are prime to $3p_j$, and so 
there is a least positive integer solution $m_j$, and this solution 
must satisfy both $m_j\leq 3p_j$ (by its minimality) and $m_j\neq 3p_j$ 
(as $(m_j, 3p_j)=1$), so that the inequalities in (3.3) will be satisfied.
\par
Since the numbers $p_1,\ldots ,p_{N+1}$ are pairwise coprime, 
the Chinese Remainder Theorem shows also that 
a solution $n\in{\Bbb Z}$ for the simultaneous congruences in (3.1) may be found, 
provided only that each one of those congruences is soluble (for $n$). 
Given any $j\in\{1,\ldots ,N+1\}$, we note that, as $(p_j , 3) = (p_j , m_j) = 1$, 
the congruence 
$3n^2\equiv m_j\pmod{p_j^2}$ is soluble (for $n$) if and only if $3m_j$ is a quadratic 
residue $\bmod\ p_j^2$. Since $p_j$ is an odd prime, this last condition on 
$3m_j$ will be satisfied if and only if $3m_j$ is a quadratic residue $\bmod\ p_j$. 
By the first congruence in (3.2), we do have $6m_j\equiv 3^2\pmod{p_j}$, so that 
$6m_j=(2)(3m_j)$ is a quadratic residue $\bmod\ p_j$. We deduce that 
the congruence $3n^2\equiv m_j\pmod{p_j^2}$ is soluble if and only if $2$ is 
a quadratic residue $\bmod\ p_j$. Given that $p_j\equiv\pm 1\pmod{8}$, 
the solubility of the congruence $3n^2\equiv m_j\pmod{p_j^2}$ therefore 
follows as a consequence of the well-known fact that,
for all odd primes $p$, the Legendre symbol $({2\over p})$ equals  
$(-1)^{(p^2 - 1)/8}$ (and so equals $1$ when $p\equiv \pm 1\pmod{8}$). 
This completes the proof of the solubility of the simultaneous 
congruences in (3.1): as $(m_j,p_j)=1$ ($j=1,\ldots,N+1$), the 
set of all integers $n$ satisfying these simultaneous congruences 
must be one of the residue classes $\bmod\ P^2$ that is prime to $P^2$, and 
so must contain a unique element $n$ that lies strictly between 
$P^2$ and $P^2 + P^2$, as in (3.4).\quad$\blacksquare$

\bigskip

\noindent{\bf The proof of Theorem~1.}\quad 
As was mentioned earlier, each eigenvalue of $K$ has a finite index 
(this is, for example, a corollary of the relations in (2.2)). 
It follows directly from this fact that $\omega^{+}$ will be some non-negative integer if 
the set ${\cal S}(K)\cap (0, \infty )$ is not infinite. 
Similarly, if the set ${\cal S}(K)\cap (-\infty , 0)$ is not infinite, then 
$\omega^{-}$ is a non-negative integer (and so not equal to $\infty$). 
Therefore Theorem~1 will follow if we can show that both 
$\omega^{-}$ and $\omega^{+}$ must equal $\infty$. We shall achieve this, in each case,  
through `proof by contradiction'. 
\par 
Suppose it is not the case that $\omega^{+}=\omega^{-}=\infty$. 
Then either $\omega^{+}\in{\Bbb N}\cup\{ 0\}$, or else $\omega^{+}=\infty$ and 
$\omega^{-}\in{\Bbb N}\cup\{ 0\}$. In the former case we put $N=\omega^{+}$ and 
$u=-1$; in the latter case we put $N=\omega^{-}$ and $u=1$. 
\par
Appealing to Lemma~3.1, we put 
$Q=5(N+1)^{1/2}$, and 
choose $N+1$ distinct primes $p_1,\ldots ,p_{N+1}>Q\geq 5$, 
with associated integers $m_1,\ldots ,m_{N+1}$ and $n$, in such a way that 
the conditions (3.1)-(3.4) are all satisfied. We then put 
$$x_j={p_j\over n}\qquad\hbox{($j=1,\ldots ,N+1$).}$$
By (3.4), the points $x_1,\ldots ,x_{N+1}$ all lie in the interval 
$(0,P/P^2)\subseteq (0,1/p_1)\subseteq (0,1/7)$.
\par 
We observe now that, as $p_j^2\equiv 1\pmod{3}$, for $j=1,\ldots,N+1$, 
the congruences (3.1) and (3.2) imply that we have $3n^2 = m_j + (3 v_j - u) p_j^2$,  
for $j=1,\ldots,N+1$, where each $v_j$ is integer valued. 
By this, we obtain:
$$K\left( x_j,x_j\right) ={1\over 2} - \left\{ {n^2\over p_j^2}\right\}
= {1\over 2} - \left\{ {m_j\over 3 p_j^2} - {u\over 3}\right\} 
= -\left( {u\over 6} + {m_j\over 3 p_j^2}\right) 
\qquad\hbox{($j=1,\ldots,N+1$),}$$
with the final equality following due to our having 
$0<m_j /(3p_j^2)<1/p_j \leq 1/7<1/3$ (as a consequence of (3.3)). 
\par 
We may note also that (3.1) and (3.2) imply that one has 
$6n^2\equiv 2m_j\equiv 3\pmod{p_j}$, and so 
$2n^2\equiv 1\pmod{p_j}$, for $j=1,\ldots, N+1$.
It follows that, when $j,k\in\{1,\ldots,N+1\}$ and $j\neq k$ (so that $(p_j,p_k)=1$), 
one will have $2n^2\equiv 1\pmod{p_j p_k}$, and so 
$2n^2 = 1 + (1+2w_{jk})p_j p_k$, where $w_{jk}$ is some integer.
We may therefore deduce that
$$K\left( x_j,x_k\right) ={1\over 2} - \left\{ {n^2\over p_j p_k}\right\}
= {1\over 2} - \left\{ {1\over 2p_j p_k} + {1\over 2}\right\} 
= - {1\over 2p_j p_k}\qquad\hbox{($1\leq j,k\leq N+1$, $j\neq k$),}$$
since we have $0<1/(2p_j p_k)<1/154<1/2$ here. 
\par 
Let now 
$$\Delta=\log\left( {t\over 1 - e^{-t}}\right)\qquad{\rm and}\qquad 
\delta = t - \Delta\;,$$
where $t>0$ is to be specified at a later point in this proof. 
By this we have 
$$e^{\Delta} - e^{-\delta} = t = \Delta + \delta\;.\eqno(3.5)$$
Note that 
$$e^t>{te^t \over e^t - 1}={t\over 1 - e^{-t}}={(t/2)e^{t/2}\over \sinh(t/2)}>{(t/2)\cosh(t/2)\over\sinh(t/2)}>1\;,$$
so that  
$$0<\Delta , \delta <t\;.$$
\par 
For $j=1,\ldots,N+1$ and $0\leq x\leq 1$, we define $\psi_j(x)$ by:  
$$\psi_j(x)=\cases{1/\sqrt{t x_j} &if $e^{-\delta}<x/x_j<e^{\Delta}$,\cr 
0 &otherwise.}$$
The functions $\psi_1(x),\ldots ,\psi_{N+1}(x)$ so defined 
are elements of the space $L^2 ([0,1])$. 
Assuming that 
$$t\leq {2\over\max\left\{ p_1,\ldots ,p_{N+1}\right\} }\;,\eqno(3.6)$$
we have 
$$\left|\log\left( {x_k\over x_j}\right)\right|
=\left|\log\left( {p_k\over p_j}\right)\right|
> {\left| p_k - p_j\right|\over\max\left\{ p_k,p_j\right\} }\geq t=\Delta + \delta
\qquad\hbox{($1\leq j,k\leq N+1$, $j\neq k$),}$$
so that, by virtue of the pairwise disjointness of the sets that are their supports,  
the functions $\psi_1,\ldots ,\psi_{N+1}$ 
are pairwise orthogonal with respect to the inner product (2.1). 
By (3.4) and (3.6), 
we have also 
$e^{\Delta}x_j < e^t p_j /n <p_j^{-1} \exp(2/p_j)\leq (1/7)\exp(2/7) < 1$, 
for $j=1,\ldots,N+1$, and so it follows (using (3.5)) that 
$\{ \psi_1,\ldots ,\psi_{N+1}\}$ is an orthonormal subset of $L^2 ([0,1])$. 
\par
Let $\psi\in L^2 ([0,1])$ be defined by 
$$\psi(x)=\sum_{j=1}^{N+1} \alpha_j \psi_j (x)\qquad\hbox{($0\leq x\leq 1$),}$$
where $\alpha_1,\ldots ,\alpha_{N+1}$ denote certain real constants that 
we shall choose later.  
Then, as a consequence of Lemma~2.1 
(combined with the fact that the square of any real number is real and non-negative),
we find that 
$$u J(\psi,\psi)=u \int_0^1 \int_0^1 K(x,y)\psi(x)\psi(y) dx dy
\geq - \sum_{j=1}^N  
{\langle \psi, \phi_j^{\pm}\rangle^2\over \left| \lambda_j^{\pm}\right|}\;,$$
where each ambiguous sign `$\pm$' is such that one has $\pm u = -1$. 
We therefore will have 
$$u J(\psi ,\psi)\geq 0\eqno(3.7)$$ 
if $\langle \psi, \phi_j^{\pm}\rangle = 0$ for $j=1,\ldots,N$. 
We observe that this last condition holds 
subject to a certain set of $N$ homogeneous linear equations in variables 
$z_1,\ldots ,z_{N+1}$ (say) being satisfied when, for $j=1,\ldots,N+1$, one has $z_j=\alpha_j$. 
The coefficients of this set of equations form an $N\times(N+1)$ real matrix, the 
columns of which are (necessarily) a linearly dependent set. 
We therefore have (3.7) for at least one choice of $\alpha_1,\ldots ,\alpha_{N+1}\in{\Bbb R}$ 
that is distinct from the `trivial solution' 
$(\alpha_1 , \ldots ,\alpha_{N+1}) = (0,\ldots ,0)$. 
We assume such a choice of $\alpha_1,\ldots,\alpha_{N+1}$ in what follows. 
Thus (3.7) holds. 
\par
By our definitions of $\psi_1,\ldots ,\psi_{N+1}$ and $\psi$, 
we find that the integral $J(\psi,\psi)$ in (3.7) satisfies
$$J(\psi,\psi) 
=\sum_{j=1}^{N+1} \sum_{k=1}^{N+1}  \Psi_{j,k}\alpha_j \alpha_k\;,$$
where
$$\Psi_{j,k} = \int_0^1 \int_0^1 K(x,y)\psi_j(x)\psi_k(y) dx dy
={1\over t \sqrt{x_j x_k}}\int_{x_j /e^{\delta}}^{x_j e^{\Delta}} 
\left( \int_{x_k /e^{\delta}}^{x_k e^{\Delta}} K(x,y) dy \right) dx\;.$$
Within the final integral here we have always
$$\left| {1\over xy} - {1\over x_j x_k}\right|<{e^{2t} -1\over x_j x_k}
={\left( e^{2t} - 1\right) n^2 \over p_j p_k} 
\leq  {\left( e^{2t} - 1\right) n^2 \over 49} \;.\eqno(3.8)$$
On the other hand, our earlier calculations of $K(x_j,x_k)$ (including that in the case $k=j$) 
make it plain that we have 
$${1\over 3} < \left\{ {1\over x_j x_k}\right\} 
= {1\over 2} - K\left( x_j , x_k\right) 
< {2\over 3} + {1\over 7} = {17\over 21}
\qquad\hbox{($1\leq j,k\leq N+1$).}$$
Therefore, provided that we choose $t>0$ so small as to satisfy 
$$e^{2t} - 1 \leq {28\over 3 n^2}\;,\eqno(3.9)$$
it will then be the case that (3.8) implies the continuity 
of the kernel $K$ at the point $(x,y)$. 
Therefore, given the particulars of the definition of $K(x,y)$, we 
may deduce (subject to (3.9)) that 
$$\Psi_{j,k} = {1\over t \sqrt{x_j x_k}}\int_{x_j /e^{\delta}}^{x_j e^{\Delta}} 
\left( \int_{x_k /e^{\delta}}^{x_k e^{\Delta}} 
\left( K\left( x_j , x_k\right) + {1\over x_j x_k} - {1\over xy}\right) dy \right) dx\;,\eqno(3.10)$$
for $1\leq j,k\leq N+1$. 
Note that (3.9) implies $t\leq 14/(3n^2)$, and so (by (3.4)), 
it is certainly a stronger condition on $t$ than that in (3.6). 
By (3.5), our result in (3.10) simplifies to:
$$\Psi_{j,k} = t \sqrt{x_j x_k} K\left( x_j,x_k\right)\qquad 
\hbox{($1\leq j,k\leq N+1$).}$$
By this, together with our earlier calculations of $K(x_j,x_k)$ (including that in the 
case $j=k$), we find that 
$$\eqalign{J(\psi,\psi) 
 &=t\sum_{j=1}^{N+1} \sum_{k=1}^{N+1}  K\left( x_j,x_k\right) \left(\alpha_j \sqrt{x_j}\right) 
\left(\alpha_k \sqrt{x_k}\right) \cr 
 &= (-t)\sum_{j=1}^{N+1} \left( {u\over 6} + {m_j\over 3 p_j^2}\right) \alpha_j^2 x_j 
 + (-t)\ \,\sum\!\sum_{\!\!\!\!\!\!\!\!\!\!\!\!\!\!\! 1\leq j < k \leq N+1} 
\left( {\alpha_j \sqrt{x_j}\over p_j}\right) 
\left( {\alpha_k \sqrt{x_k}\over p_k}\right)\;, \cr }$$
and so 
$$u J(\psi , \psi) 
= - {t\over n}\left( 
\sum_{j=1}^{N+1} \left( {p_j\over 6} + {u m_j\over 3 p_j}\right) \alpha_j^2 
 + u\,\sum\!\sum_{\!\!\!\!\!\!\!\!\!\!\!\!\!\!\! 1\leq j < k \leq N+1} 
{\alpha_j \alpha_k \over \sqrt{p_j p_k}}
\right) \;.$$
Here the absolute value of the sum over $j$ and $k$ may be bounded 
by applying the triangle inequality and the inequality 
$\left| \alpha_j\alpha_k /\sqrt{p_j p_k}\right|\leq 
\left( \alpha_j^2 p_j^{-1} + \alpha_k^2 p_k^{-1}\right)/2$: 
in this way one obtains, with the help of (3.3),  the upper bound  
$$u J(\psi,\psi)\leq 
- {t\over n} 
\sum_{j=1}^{N+1} \left( {p_j\over 6} -\left( 1 + {N\over 2 p_j}\right)\right) \alpha_j^2 \;. 
\eqno(3.11)$$
Since we have $\min\{ p_1,\ldots ,p_{N+1}\} > Q\geq 5\max\{ 1 , \sqrt{N}\}$, 
if follows that 
$$1/\left( {p_j\over 6}\right) ={6\over p_j}\leq {6\over 7}\quad 
{\rm and}\quad 
\left( {N\over 2 p_j}\right) / \left( {p_j\over 6}\right)  = {3N\over p_j^2} 
< {3N\over Q^2} \leq {3\over 25} < {1\over 8}\;,$$ 
for $j=1,\ldots , N+1$, and so, by (3.11) and the `non-triviality' of 
$(\alpha_1 , \ldots ,\alpha_{N+1})\in{\Bbb R}^{N+1}$, 
$$u J(\psi , \psi) 
\leq - {t\over n} 
\sum_{j=1}^{N+1}  {p_j \alpha_j^2\over 336} 
\leq  - {t\over 48 n} \sum_{j=1}^{N+1} \alpha_j^2 < 0\;.
$$
Since the last of these inequalities is strict, we find that (3.7) is contradicted, and so 
complete our `proof by contradiction' that 
$\omega^{+}=\omega^{-}=\infty$.$\quad\blacksquare$
\medskip 
\noindent{\bf Remarks.} 

\smallskip

\noindent{\bf (1)}\ The idea for the above proof came after reading some of H. Weyl's paper [5]: 
in particular, his proof of `Satz~1' there. 

\smallskip

\noindent{\bf (2)}\ By elaborating upon the above proof 
one can obtain lower bounds for the terms in the sequence $(\lambda_k^{-})$, 
and upper bounds for the terms in the sequence $(\lambda_k^{+})$. 
These bounds, however, are extremely weak: the best I have been able to show, 
in respect of positive eigenvalues of $K$, is that one has 
$$\lambda_k^{+}\leq 2772 \bigl( 918(k+1)\log(k+1)\bigr)^{6(k+1)}\qquad 
\hbox{($k\in{\Bbb N}$).}$$

\bigskip

\beginsection References.  

\item{[1]} M.N. Huxley \& N. Watt, ``Mertens sums requiring fewer values of the M\"obius function'', 
{\it Chebyshevski\u{i} Sb.} {\bf 19} (2018), no. 3, p. 19-34. 

\item{[2]} F. Mertens, ``\"{U}ber eine zahlentheoretische Function'', 
{\it Sitzungsber. Ak. Wiss. Wien, Math.-Naturw. Kl., Abt. IIa} {\bf 106} (1897), p. 761-830. 

\item{[3]} F.G. Tricomi, {\it Integral Equations}, Dover Publications, Inc., 1985, viii+238 pages. 
Reprint. Originally Published: Interscience Publishers, 1957.  

\item{[4]} N. Watt, ``The kernel ${1\over 2} +\lfloor {1\over xy}\rfloor - {1\over xy}\,$ ($0<x,y\leq 1$) 
and Mertens sums'', \hfill{\tt https://arxiv}\break{\tt .org/abs/1812.01039v1}, 2018. 

\item{[5]} H. Weyl, ``Ueber die asymptotische Verteilung der Eigenwerte'', 
{\it Nachr. Ges. Wiss. G\"{o}ttingen, Math.-Phys. Kl.} (1911), p. 110-117.

\vskip 1cm

\bye